\documentclass{elsart}
\usepackage{graphics}
\usepackage{graphicx}
\usepackage{epsfig}
\usepackage{amssymb,amsmath}
\newtheorem{Lemma}{Lemma}[section]
\newtheorem{Theorem}{Theorem}[section]
\newtheorem{Remark}{Remark}[section]
\begin{document}
\begin{frontmatter}
\title{Lyapunov-Sylvester Operators for Numerical Solutions of NLS Equation}
\author{Riadh CHTEOUI}
\ead{riadh.chteoui2005@yahoo.fr}
\address{D\'epartement de Math\'ematiques, Facult\'e des Sciences de Monastir, 5019 Monastir, Tunisia.}
\author{Anouar BEN MABROUK\corauthref{cor1}\thanksref{label2}}
\ead{anouar.benmabrouk@issatso.rnu.tn}
\corauth[cor1]{Corresponding Author}
\address{Computational Mathematics Lab, D\'epartement de Math\'ematiques, Facult\'e des Sciences, 5019 Monastir, Tunisia.}
\thanks[label2]{D\'epartement de Math\'ematiques, Institut Sup\'erieur de Math\'amatiques Appliqu\'ees et Informatique de Kairouan, Avenue Assad Ibn Al-Fourat, Kairouan 3100, Tunisia.}
\author{Hichem OUNAIES}
\ead{hichem.ounaies@fsm.rnu.tn}
\address{D\'epartement de Math\'ematiques, Facult\'e des Sciences de Monastir, 5019 Monastir, Tunisia.}
\begin{abstract}
In the present paper a numerical method is developed to approximate the solution of two-dimensional NLS equation in the presence of a singular potential. The method leads to Lyapunov-Syslvester algebraic operators that are shown to be invertible using original topological and differential calculus issued methods. The numerical scheme is proved to be consistent, convergent and stable using the based on Lyapunov criterion, lax equivalence theorem and the properties of the Lyapunov-Syslvester operators.
\end{abstract}
\end{frontmatter}
\section{Introduction}
The Schr\"odinger equation is widely studied from both numerical and theoretical points of view. This is due to its relation to the modeling of real physical phenomena such as Newton's laws and conservation of energy in classical mechanics, behaviour of dynamical systems, the description of a particle in a non-relativistic setting in quantum mechanics, etc. The Schr\"odinger's linear equation states that
$$\Delta\psi +\frac{8\pi^2m}{\hbar^2}\,\left(E-V(x)\right)\psi =0\,,$$
where $\psi$ is the Schr\"odinger wave function, $m$ is the mass, $\hbar$ denotes Planck's constant, $E$ is the energy, and $V$ is the potential energy. However, in the nonlinear case, the structure of the nonlinear Schr\"odinger equation is more complicated. This equation is a prototypical dispersive nonlinear partial differential equation related Bose-Einstein condensates and nonlinear optics (\cite{byeon}), propagation of electric fields in optical fibers (\cite{hase}, \cite{malo}), self-focusing and collapse of Langmuir waves in plasma physics (\cite{zakh}), behaviour of rogue waves in oceans (\cite{onor}).

The nonlinear Schr\"odinger equation is also related to electromagnetic, ferromagnetic fields as well as magnums, high-power ultra-short laser self-channelling in matter, condensed matter theory, dissipative quantum mechanics, (\cite{avron1}), film equations, etc (See \cite{abl}, \cite{sulem}).

Based upon the analogy between mechanics and optics, Schr\"odinger established the classical derivation of his equation. By developing a perturbation method, he proved the equivalence between his wave mechanics equation and and Heisenberg's matrix one, and thus introduced the time dependent version stated hereafter with a cubic nonlinearity
\begin{equation}\label{oheq}
i\hbar\psi_t=-\frac{\hbar^2}{2m}\,\Delta\psi+V(x)\psi-\gamma |\psi|^{2}\psi\qquad\mbox{in $\mathbb{R}^N$ ($N\geq 2$)}.
\end{equation}
In \cite{floer} and \cite{oh} the potential $V$ is assumed to be bounded with a non-degenerate critical point at $x=0$. More precisely, $V$ belongs to the class $V_a$, for some real parameter $a$ (See \cite{kato}). With suitable assumptions it is proved in \cite{oh} a Lyapunov-Schmidt type reduction the existence of standing wave solutions of problem \eqref{oheq}, of the form
\begin{equation}\label{oheq1}
\psi (x,t)=e^{-iEt/\hbar}u(x)\,.
\end{equation}
Hence, the nonlinear Schr\"odinger equation \eqref{oheq} is reduced to the semilinear elliptic equation
$$
-\frac{\hbar^2}{2m}\,\Delta u+\left(V(x)-E\right)u=|u|^{2}u\,.
$$
Setting $y=\hbar^{-1}x$ and replacing $y$ by $x$ we get
\begin{equation}\label{oheqoheq}
-\Delta u+2m\left(V_\hbar (x)-E\right)u=|u|^{2}u\qquad\mbox{in $\mathbb{R}^N$}\,,
\end{equation}
where $V_\hbar (x)=V(\hbar x)$.

If for some $\xi\in\mathbb{R}^N\setminus\{0\}$, $V(x + s\xi) = V (x)$ for all $s \in\mathbb{R}$, equation \eqref{oheq} is invariant under the Galilean transformation
$$
\psi(x,t)\longmapsto\exp\left(i\xi\cdot x/\hbar -\frac12i|\xi|^2t/\hbar\right)\psi(x-\xi t, t)\,.
$$
Thus, in this case, standing waves reproduce solitary waves traveling in the direction of $\xi$.

The present paper is devoted to the development of a numerical method based on two-dimensional finite difference scheme to approximate the solution of the nonlinear Schr\"odinger (NLS) equation in $\mathbb{R}^2$ written on the form
\begin{equation}\label{eqn1-1}
iu_{t}+\Delta\,u+u-|u|^{-2\theta}u=0,\quad\in\Omega\times(t_0,+\infty)
\end{equation}
with the initial and boundary conditions
\begin{equation}\label{eqn1-2}
u(x,y,t_0)=u_0(x,y)\quad\hbox{and}\quad\displaystyle\frac{\partial\,u}{\partial\,n}(x,y,t)=0,\quad(x,y,t)\in\Omega\times(t_0,+\infty).
\end{equation}
We consider a rectangular domain $\Omega=]L_0,L_1[\times]L_0,L_1[$ in $\mathbb{R}^2$ and $t_0$ a real parameter fixed as the initial time, $u_{t}$ is the first order partial derivative in time, $\Delta=\frac{\partial^2}{\partial\,x^2}+\frac{\partial^2}{\partial\,y^2}$ is the Laplace operator in $\mathbb{R}^2$.
$\displaystyle\frac{\partial}{\partial n}$ is the outward normal derivative operator along the boundary $\partial\Omega$. Finally,
$u$ and $u_0$ are complexe valued functions.

In \cite{chteouietal}, the stationary solutions of problem (\ref{eqn1-1}) has been studied using direct methods issued from the equation on the whole space. Existence, uniqueness, classification and properties of the solutions have been investigated. It is proved that three attractive zones or three classes of stationary solutions exists; there are solutions oscillating around 0 with supports being compact, there are solutions oscillating around $\pm1$ with a finite number of zeros.
In the present paper, we reconsider the evolutionary NLS equation associated to the stationary problem studied in \cite{chteouietal}. The idea consists in developing a numerical scheme to approximate the solution of (\ref{eqn1-1})-(\ref{eqn1-2}). The time and the space partial derivatives are replaced by barycenter finite-difference approximations in order to transform the initial boundary-value problem (\ref{eqn1-1})-(\ref{eqn1-2}) into a linear algebraic system. The resulting method is analyzed for local truncation error and stability and the scheme is proved to be uniquely solvable and convergent.

We recall that in the majority of previous works Lyapunov type equations are not widely applied but standard linear tridiagonal and/or fringe-tridiagonal operators acting on one column vector. This leads to the application of classical methods such that Shur decomposition, relaxation, etc to compute the inverse of linear operators obtained. In \cite{Benmabrouk-Ayadi2} a first investigation was done based on Lyapunov operators to solve numerical NLS and Heat equations. It is proved that the Lyapunov operators obtained may be transformed with translation-dilation actions into contractive ones which leads to uniquely solvable systems using fixed point theory. Here. we do not apply the same computations as in \cite{Benmabrouk-Ayadi2}, but we develop different arguments. We will instead apply a differential calculus and topology technique (See \cite{HenriCartan} for example) to prove theorem \ref{theorem1}. We recall finally that constants appearing along the proofs are generic and may differ from one step to an other and are denoted undifferently by the capital $C$.
\section{Discrete Two-Dimensional NLS Equation}
Let $\Omega=]L_0,L_1[\times]L_0,L_1[\subset\mathbb{R}^2$ and for $J\in\mathbb{N}^*$, denote $h=\displaystyle\frac{L_1-L_0}{J}$ for the space step, $x_j=L_0+jh$ and $y_m=L_0+mh$ for all $(j,m)\in{I}^2=\{0,1,\dots,J\}^2$. Let $l=\Delta\,t$ be the time step and $t_n=t_0+nl$, $n\in\mathbb{N}$ be the discrete time grid. For $(j,m)\in{I}$ and $n\geq0$, $u_{j,m}^n$ will be the net function $u(x_j,y_m,t_n)$ and $U_{j,m}^n$ the numerical solution. The following discrete approximations will be applied for the different differential operators involved in the problem. For time derivatives, we set
$$
\displaystyle\,u_t\rightsquigarrow\displaystyle\frac{U_{j,m}^{n+1}-U_{j,m}^{n-1}}{2l}
$$
and for space derivatives, we shall use
$$
\displaystyle\,u_x\rightsquigarrow\displaystyle\frac{U_{j+1,m}^{n}-U_{j-1,m}^{n}}{2h}\quad\hbox{and}\quad
\displaystyle\,u_y\rightsquigarrow\displaystyle\frac{U_{j,m+1}^{n}-U_{j,m-1}^{n}}{2h}
$$
for first order derivatives and for second order ones we apply the estimations
$$
\Delta\,u\rightsquigarrow\Delta\,U^n_{j,m}(\alpha_n,\beta_n,\gamma_n)=\Delta\biggl[\alpha_nU_{j,m}^{n+1}+\beta_nU_{j,m}^{n}+\gamma_nU_{j,m}^{n-1}\biggr]
$$
where
$$
\Delta\,U_{j,m}^{n}=\displaystyle\frac{U_{j+1,m}^{n}-2U_{j,m}^{n}+U_{j-1,m}^{n}}{h^2}+\displaystyle\frac{U_{j,m+1}^{n}-2U_{j,m}^{n}+U_{j,m-1}^{n}}{h^2}.
$$
where $\alpha_n$, $\beta_n$ et $\gamma_n$ are sequences in [0,1] such that $\alpha_n+\beta_n+\gamma_n=1$.

By replacing the derivatives of $u$ with their approximations, the equation (\ref{eqn1-1}) yields
$$
i\displaystyle\frac{U_{j,m}^{n+1}-U_{j,m}^{n-1}}{2l}+\Delta\biggl[\alpha_nU_{j,m}^{n+1}+\beta_nU_{j,m}^{n}+\gamma_nU_{j,m}^{n-1}\biggr]+\widehat{f_{j,m}^{n}}=0
$$
where $\widehat{f_{j,m}^{n}}=\displaystyle\frac{f(U_{j,m}^{n})+f(U_{j,m}^{n-1})}{2}$ where we design by $f(u)=u-|u|^{-2\theta}u$. We then obtain
$$
\begin{matrix}\medskip\,i\displaystyle\frac{U_{j,m}^{n+1}-U_{j,m}^{n-1}}{2l}+\alpha_n\Biggl[\displaystyle\frac{U_{j-1,m}^{n+1}-2U_{j,m}^{n+1}+U_{j+1,m}^{n+1}
+U_{j,m-1}^{n+1}-2U_{j,m}^{n+1}+U_{j,m+1}^{n+1}}{h^2}\Biggr]\hfill\cr\medskip
+\beta_n\Biggl[\displaystyle\frac{U_{j-1,m}^{n}-2U_{j,m}^{n}+U_{j+1,m}^{n}+U_{j,m-1}^{n}-2U_{j,m}^{n}+U_{j,m+1}^{n}}{h^2}\Biggr]\hfill\cr\medskip
+\gamma_n\Biggl[\displaystyle\frac{U_{j-1,m}^{n-1}-2U_{j,m}^{n-1}+U_{j+1,m}^{n-1}+U_{j,m-1}^{n-1}-2U_{j,m}^{n-1}+U_{j,m+1}^{n-1}}{h^2}\Biggr]
+\widehat{f_{j,m}^{n}}=0.\hfill\end{matrix}
$$
Denote next $\sigma=\displaystyle\frac{2l}{h^2}$. We obtain
$$
\begin{matrix}\medskip\,iU_{j,m}^{n+1}-iU_{j,m}^{n-1}+\sigma\alpha_n\Biggl[U_{j-1,m}^{n+1}-2U_{j,m}^{n+1}+U_{j+1,m}^{n+1}
+U_{j,m-1}^{n+1}-2U_{j,m}^{n+1}+U_{j,m+1}^{n+1}\Biggr]\hfill\cr\medskip
+\sigma\beta_n\Biggl[U_{j-1,m}^{n}-2U_{j,m}^{n}+U_{j+1,m}^{n}+U_{j,m-1}^{n}-2U_{j,m}^{n}+U_{j,m+1}^{n}\Biggr]\hfill\cr\medskip
+\sigma\gamma_n\Biggl[U_{j-1,m}^{n-1}-2U_{j,m}^{n-1}+U_{j+1,m}^{n-1}+U_{j,m-1}^{n-1}-2U_{j,m}^{n-1}+U_{j,m+1}^{n-1}\Biggr]
+2l\widehat{f_{j,m}^{n}}=0.\hfill\end{matrix}
$$
By setting $\varphi_n=\displaystyle\frac{i-4\sigma\alpha_n}{2}$ and $\psi_n=\displaystyle\frac{i+4\sigma\gamma_n}{2}$ and gathering the terms according to the time steps we obtain
$$
\begin{matrix}\medskip\,
\sigma\alpha_nU_{j-1,m}^{n+1}+\varphi_nU_{j,m}^{n+1}+\sigma\alpha_nU_{j+1,m}^{n+1}\hfill\cr\medskip
+\sigma\alpha_n U_{j,m-1}^{n+1}+\varphi_nU_{j,m}^{n+1}+\sigma\alpha_nU_{j,m+1}^{n+1}\hfill\cr\medskip
+\sigma\beta_nU_{j-1,m}^{n}-2\sigma\beta_nU_{j,m}^{n}+\sigma\beta_nU_{j+1,m}^{n}\hfill\cr\medskip
+\sigma\beta_nU_{j,m-1}^{n}-2\sigma\beta_nU_{j,m}^{n}+\sigma\beta_nU_{j,m+1}^{n}\hfill\cr\medskip
+\sigma\alpha_nU_{j-1,m}^{n-1}-\psi_nU_{j,m}^{n-1}+\sigma\alpha_nU_{j+1,m}^{n-1}\hfill\cr\medskip
+\sigma\alpha_nU_{j,m-1}^{n-1}-\psi_nU_{j,m}^{n-1}+\sigma\alpha_nU_{j,m+1}^{n-1}+2l\widehat{f_{j,m}^{n}}=0.\hfill\end{matrix}
$$
or in vector form
$$
\begin{matrix}
\left(\begin{array}{ccc}
\sigma\alpha_n & \varphi_n & \sigma\alpha_n \\
\end{array}\right)
\left(
  \begin{array}{c}
    U_{j-1,m}^{n+1} \\
    U_{j,m}^{n+1} \\
    U_{j+1,m}^{n+1} \\
  \end{array}
\right)
+
\left(\begin{array}{ccc}
U_{j,m-1}^{n+1} & U_{j,m}^{n+1} & U_{j,m+1}^{n+1} \\
\end{array}\right)
\left(
  \begin{array}{c}
    \sigma\alpha_n \\
    \varphi_n \\
    \sigma\alpha_n \\
  \end{array}
\right)\hfill\cr
+\sigma\beta_n
\left(\begin{array}{ccc}
1 & -2 & 1 \\
\end{array}\right)
\left(
  \begin{array}{c}
    U_{j-1,m}^{n} \\
    U_{j,m}^{n} \\
    U_{j+1,m}^{n} \\
  \end{array}
\right)
+\sigma\beta_n
\left(\begin{array}{ccc}
U_{j,m-1}^{n} & U_{j,m}^{n} & U_{j,m+1}^{n} \\
\end{array}\right)
\left(
  \begin{array}{c}
    1 \\
    -2 \\
    1 \\
  \end{array}
\right)\hfill\cr
+\left(\begin{array}{ccc}
\sigma\gamma_n & \psi_n & \sigma\gamma_n \\
\end{array}\right)
\left(
  \begin{array}{c}
    U_{j-1,m}^{n-1} \\
    U_{j,m}^{n-1} \\
    U_{j+1,m}^{n-1} \\
  \end{array}
\right)
+
\left(\begin{array}{ccc}
U_{j,m-1}^{n-1} & U_{j,m}^{n-1} & U_{j,m+1}^{n-1} \\
\end{array}\right)
\left(
  \begin{array}{c}
    \sigma\gamma_n \\
    \psi_n \\
    \sigma\gamma_n \\
  \end{array}
\right)+2l\widehat{f_{j,m}^{n}}=0.\hfill
\end{matrix}
$$
Now, we exploit the boundary conditions which can be resumed in the following cases of the parameters $j,m$. Indeed, by setting in the previous equation $j=m=0$ and using the approximations of boundary conditions we obtain
$$
\begin{matrix}
\left(\begin{array}{cc}
\varphi_n & 2\sigma\alpha_n \\
\end{array}\right)
\left(
  \begin{array}{c}
    U_{0,0}^{n+1} \\
    U_{1,0}^{n+1} \\
  \end{array}
\right)
+
\left(\begin{array}{cc}
U_{0,0}^{n+1} & U_{0,1}^{n+1} \\
\end{array}\right)
\left(
  \begin{array}{c}
    \varphi_n \\
    2\sigma\alpha_n \\
  \end{array}
\right)\hfill\cr
+\sigma\beta_n
\left(\begin{array}{cc}
-2 & 2 \\
\end{array}\right)
\left(
  \begin{array}{c}
    U_{0,0}^{n} \\
    U_{1,0}^{n} \\
  \end{array}
\right)
+\sigma\beta_n
\left(\begin{array}{cc}
U_{0,0}^{n} & U_{0,1}^{n} \\
\end{array}\right)
\left(
  \begin{array}{c}
    -2 \\
    2 \\
  \end{array}
\right)\hfill\cr
+\left(\begin{array}{cc}
\psi_n & 2\sigma\gamma_n \\
\end{array}\right)
\left(
  \begin{array}{c}
    U_{0,0}^{n-1} \\
    U_{1,0}^{n-1} \\
  \end{array}
\right)
+
\left(\begin{array}{cc}
U_{0,0}^{n-1} & U_{0,1}^{n-1} \\
\end{array}\right)
\left(
  \begin{array}{c}
    \psi_n \\
    2\sigma\gamma_n \\
  \end{array}
\right)+2l\widehat{f_{0,0}^{n}}=0.\hfill
\end{matrix}
$$
For $j=0$ and $m=J$, we obtain as previously,
$$
\begin{matrix}
\left(\begin{array}{cc}
\varphi_n & 2\sigma\alpha_n \\
\end{array}\right)
\left(
  \begin{array}{c}
    U_{0,J}^{n+1} \\
    U_{1,J}^{n+1} \\
  \end{array}
\right)
+
\left(\begin{array}{cc}
U_{0,J-1}^{n+1} & U_{0,J}^{n+1} \\
\end{array}\right)
\left(
  \begin{array}{c}
    2\sigma\alpha_n \\
    \varphi_n \\
  \end{array}
\right)\hfill\cr
+\sigma\beta_n
\left(\begin{array}{cc}
-2 & 2 \\
\end{array}\right)
\left(
  \begin{array}{c}
    U_{0,J}^{n} \\
    U_{1,J}^{n} \\
  \end{array}
\right)
+\sigma\beta_n
\left(\begin{array}{cc}
U_{0,J-1}^{n} & U_{0,J}^{n} \\
\end{array}\right)
\left(
  \begin{array}{c}
    2 \\
    -2 \\
  \end{array}
\right)\hfill\cr
+\left(\begin{array}{cc}
\psi_n & 2\sigma\gamma_n \\
\end{array}\right)
\left(
  \begin{array}{c}
    U_{0,J}^{n-1} \\
    U_{1,J}^{n-1} \\
  \end{array}
\right)
+
\left(\begin{array}{cc}
U_{0,J-1}^{n-1} & U_{0,J}^{n-1} \\
\end{array}\right)
\left(
  \begin{array}{c}
     2\sigma\gamma_n \\
     \psi_n \\
  \end{array}
\right)+2l\widehat{f_{0,J}^{n}}=0.\hfill
\end{matrix}
$$
For $j=0$ and $1\leq m\leq J-1$, we get
$$
\begin{matrix}
\left(\begin{array}{cc}
\varphi_n & 2\sigma\alpha_n \\
\end{array}\right)
\left(
  \begin{array}{c}
    U_{0,m}^{n+1} \\
    U_{1,m}^{n+1} \\
  \end{array}
\right)
+
\left(\begin{array}{ccc}
U_{0,m-1}^{n+1} & U_{0,m}^{n+1} & U_{0,m+1}^{n+1} \\
\end{array}\right)
\left(
  \begin{array}{c}
    \sigma\alpha_n \\
    \varphi_n \\
    \sigma\alpha_n \\
  \end{array}
\right)\hfill\cr
+\sigma\beta_n
\left(\begin{array}{cc}
-2 & 2 \\
\end{array}\right)
\left(
  \begin{array}{c}
    U_{0,m}^{n} \\
    U_{1,m}^{n} \\
  \end{array}
\right)
+\sigma\beta_n
\left(\begin{array}{ccc}
U_{0,m-1}^{n} & U_{0,m}^{n} & U_{0,m+1}^{n} \\
\end{array}\right)
\left(
  \begin{array}{c}
    1 \\
    -2 \\
    1 \\
  \end{array}
\right)\hfill\cr
+\left(\begin{array}{cc}
\psi_n & 2\sigma\gamma_n \\
\end{array}\right)
\left(
  \begin{array}{c}
    U_{0,m}^{n-1} \\
    U_{1,m}^{n-1} \\
  \end{array}
\right)
+
\left(\begin{array}{ccc}
U_{0,m-1}^{n-1} & U_{0,m}^{n-1} & U_{0,m+1}^{n-1}\\
\end{array}\right)
\left(
  \begin{array}{c}
    \sigma\gamma_n \\
    \psi_n \\
    \sigma\gamma_n \\
  \end{array}
\right)+2l\widehat{f_{0,m}^{n}}=0.\hfill
\end{matrix}
$$
Now, by setting $j=J$ and $m=0$ and using the approximations of boundary conditions we obtain
$$
\begin{matrix}
\left(\begin{array}{cc}
2\sigma\alpha_n & \varphi_n \\
\end{array}\right)
\left(
  \begin{array}{c}
    U_{J-1,0}^{n+1} \\
    U_{J,0}^{n+1} \\
  \end{array}
\right)
+
\left(\begin{array}{cc}
U_{J,0}^{n+1} & U_{J,1}^{n+1} \\
\end{array}\right)
\left(
  \begin{array}{c}
    \varphi_n \\
    2\sigma\alpha_n \\
  \end{array}
\right)\hfill\cr
+\sigma\beta_n
\left(\begin{array}{cc}
2 & -2 \\
\end{array}\right)
\left(
  \begin{array}{c}
    U_{J-1,0}^{n} \\
    U_{J,0}^{n} \\
  \end{array}
\right)
+\sigma\beta_n
\left(\begin{array}{cc}
U_{J,0}^{n} & U_{J,1}^{n} \\
\end{array}\right)
\left(
  \begin{array}{c}
    -2 \\
    2 \\
  \end{array}
\right)\hfill\cr
+\left(\begin{array}{cc}
2\sigma\gamma_n & \psi_n \\
\end{array}\right)
\left(
  \begin{array}{c}
    U_{J-1,0}^{n-1} \\
    U_{J,0}^{n-1} \\
  \end{array}
\right)
+
\left(\begin{array}{cc}
U_{J,0}^{n-1} & U_{J,1}^{n-1} \\
\end{array}\right)
\left(
  \begin{array}{c}
    \psi_n \\
    2\sigma\gamma_n \\
  \end{array}
\right)+2l\widehat{f_{J,0}^{n}}=0.\hfill
\end{matrix}
$$
For $j=J$ and $m=J$, we obtain
$$
\begin{matrix}
\left(\begin{array}{cc}
2\sigma\alpha_n & \varphi_n \\
\end{array}\right)
\left(
  \begin{array}{c}
    U_{J-1,J}^{n+1} \\
    U_{J,J}^{n+1} \\
  \end{array}
\right)
+
\left(\begin{array}{cc}
U_{J,J-1}^{n+1} & U_{J,J}^{n+1} \\
\end{array}\right)
\left(
  \begin{array}{c}
    2\sigma\alpha_n \\
    \varphi_n \\
  \end{array}
\right)\hfill\cr
+\sigma\beta_n
\left(\begin{array}{cc}
2 & -2 \\
\end{array}\right)
\left(
  \begin{array}{c}
    U_{J-1,J}^{n} \\
    U_{J,J}^{n} \\
  \end{array}
\right)
+\sigma\beta_n
\left(\begin{array}{cc}
U_{J,J-1}^{n} & U_{J,J}^{n} \\
\end{array}\right)
\left(
  \begin{array}{c}
    2 \\
    -2 \\
  \end{array}
\right)\hfill\cr
+\left(\begin{array}{cc}
2\sigma\gamma_n & \psi_n \\
\end{array}\right)
\left(
  \begin{array}{c}
    U_{J-1,J}^{n-1} \\
    U_{J,J}^{n-1} \\
  \end{array}
\right)
+
\left(\begin{array}{cc}
U_{J,J-1}^{n-1} & U_{J,J}^{n-1} \\
\end{array}\right)
\left(
  \begin{array}{c}
     2\sigma\gamma_n \\
     \psi_n \\
  \end{array}
\right)+2l\widehat{f_{J,J}^{n}}=0.\hfill
\end{matrix}
$$
For $j=J$ and $1\leq m\leq J-1$ there holds that
$$
\begin{matrix}
\left(\begin{array}{cc}
2\sigma\alpha_n & \varphi_n \\
\end{array}\right)
\left(
  \begin{array}{c}
    U_{J-1,m}^{n+1} \\
    U_{J,m}^{n+1} \\
  \end{array}
\right)
+
\left(\begin{array}{ccc}
U_{J,m-1}^{n+1} & U_{J,m}^{n+1} & U_{J,m+1}^{n+1} \\
\end{array}\right)
\left(
  \begin{array}{c}
    \sigma\alpha_n \\
    \varphi_n \\
    \sigma\alpha_n \\
  \end{array}
\right)\hfill\cr
+\sigma\beta_n
\left(\begin{array}{cc}
2 & -2 \\
\end{array}\right)
\left(
  \begin{array}{c}
    U_{J-1,m}^{n} \\
    U_{J,m}^{n} \\
  \end{array}
\right)
+\sigma\beta_n
\left(\begin{array}{ccc}
U_{J,m-1}^{n} & U_{J,m}^{n} & U_{J,m+1}^{n} \\
\end{array}\right)
\left(
  \begin{array}{c}
    1 \\
    -2 \\
    1 \\
  \end{array}
\right)\hfill\cr
+\left(\begin{array}{cc}
2\sigma\gamma_n & \psi_n \\
\end{array}\right)
\left(
  \begin{array}{c}
    U_{J-1,m}^{n-1} \\
    U_{J,m}^{n-1} \\
  \end{array}
\right)
+
\left(\begin{array}{ccc}
U_{J,m-1}^{n-1} & U_{J,m}^{n-1} & U_{J,m+1}^{n-1}\\
\end{array}\right)
\left(
  \begin{array}{c}
    \sigma\gamma_n \\
    \psi_n \\
    \sigma\gamma_n \\
  \end{array}
\right)+2l\widehat{f_{J,m}^{n}}=0.\hfill
\end{matrix}
$$
Next, for $1\leq j\leq J-1$ and $m=0$, we have
$$
\begin{matrix}
\left(\begin{array}{ccc}
\sigma\alpha_n & \varphi_n & \sigma\alpha_n\\
\end{array}\right)
\left(
  \begin{array}{c}
    U_{j-1,0}^{n+1} \\
    U_{j,0}^{n+1} \\
    U_{j+1,0}^{n+1} \\
  \end{array}
\right)
+
\left(\begin{array}{cc}
U_{j,0}^{n+1} & U_{j,1}^{n+1} \\
\end{array}\right)
\left(
  \begin{array}{c}
    \varphi_n \\
    2\sigma\alpha_n \\
  \end{array}
\right)\hfill\cr
+\sigma\beta_n
\left(\begin{array}{ccc}
1& -2 & 1 \\
\end{array}\right)
\left(
  \begin{array}{c}
    U_{j-1,0}^{n} \\
    U_{j,0}^{n} \\
    U_{j+1,0}^{n} \\
  \end{array}
\right)
+\sigma\beta_n
\left(\begin{array}{cc}
U_{j,0}^{n} & U_{j,1}^{n} \\
\end{array}\right)
\left(
  \begin{array}{c}
    -2 \\
    2 \\
  \end{array}
\right)\hfill\cr
+\left(\begin{array}{ccc}
\sigma\gamma_n & \psi_n & \sigma\gamma_n \\
\end{array}\right)
\left(
  \begin{array}{c}
    U_{j-1,0}^{n-1} \\
    U_{j,0}^{n-1} \\
    U_{j+1,0}^{n-1} \\
  \end{array}
\right)
+
\left(\begin{array}{cc}
U_{j,0}^{n-1} & U_{j,1}^{n-1} \\
\end{array}\right)
\left(
  \begin{array}{c}
    \psi_n \\
    2\sigma\gamma_n \\
  \end{array}
\right)+2l\widehat{f_{j,0}^{n}}=0.\hfill
\end{matrix}
$$
Now, for $1\leq j\leq J-1$ and $m=J$, we obtain
$$
\begin{matrix}
\left(\begin{array}{ccc}
\sigma\alpha_n & \varphi_n & \sigma\alpha_n\\
\end{array}\right)
\left(
  \begin{array}{c}
    U_{j-1,J}^{n+1} \\
    U_{j,J}^{n+1} \\
    U_{j+1,J}^{n+1} \\
  \end{array}
\right)
+
\left(\begin{array}{cc}
U_{j,J-1}^{n+1} & U_{j,J}^{n+1} \\
\end{array}\right)
\left(
  \begin{array}{c}
    2\sigma\alpha_n \\
    \varphi_n \\
  \end{array}
\right)\hfill\cr
+\sigma\beta_n
\left(\begin{array}{ccc}
1 & -2 & 1 \\
\end{array}\right)
\left(
  \begin{array}{c}
    U_{j-1,J}^{n} \\
    U_{j,J}^{n} \\
    U_{j+,J}^{n} \\
  \end{array}
\right)
+\sigma\beta_n
\left(\begin{array}{cc}
U_{j,J-1}^{n} & U_{J,J}^{n} \\
\end{array}\right)
\left(
  \begin{array}{c}
    2 \\
    -2 \\
  \end{array}
\right)\hfill\cr
+\left(\begin{array}{ccc}
\sigma\gamma_n & \psi_n & \sigma\gamma_n\\
\end{array}\right)
\left(
  \begin{array}{c}
    U_{j-1,J}^{n-1} \\
    U_{j,J}^{n-1} \\
    U_{j+1,J}^{n-1} \\
  \end{array}
\right)
+
\left(\begin{array}{cc}
U_{j,J-1}^{n-1} & U_{j,J}^{n-1} \\
\end{array}\right)
\left(
  \begin{array}{c}
     2\sigma\gamma_n \\
     \psi_n \\
  \end{array}
\right)+2l\widehat{f_{j,J}^{n}}=0.\hfill
\end{matrix}
$$
In the matrix form we obtain
\begin{equation}\label{Lyap1}
    \mathcal{L}_{A_n}(U^{n+1})+\mathcal{L}_{B_n}(U^{n})+\mathcal{L}_{C_n}(U^{n-1})+F_n=0,
\end{equation}
where $U^n=\left(U_{j,m}^n\right)$ is the matrix of the numerical solution. $A_n$, $B_n$ and $C_n$ are $(J+1)\times(J+1)$ tri-diagonal matrices with respective coefficients
$$
\begin{array}{c}
A_n(j,j)=\varphi_n,\\
A_n(j,j-1)=A_n(j,j+1)=\sigma\alpha_n,\\
A_n(0,1)=A_n(J,J-1)=2\sigma\alpha_n.\\
\end{array}
$$
$$
\begin{array}{c}
B_n(j,j)=-2\sigma\beta_n,\\
B_n(j,j-1)=B_n(j,j+1)=\sigma\beta_n,\\
B_n(0,1)=B_n(J,J-1)=2\sigma\beta_n.\\
\end{array}
$$
$$
\begin{array}{c}
C_n(j,j)=\psi_n,\\
C_n(j,j-1)=C_n(j,j+1)=\sigma\gamma_n,\\
C_n(0,1)=C_n(J,J-1)=2\sigma\gamma_n.\\
\end{array}
$$
and $F^n=(2l\widehat{f_{j,m}^n})$. Finally, for a matrix $W\in\mathcal{M}_{J+1}(\mathbb{C})$, $\mathcal{L}_W$ is the so-called Lyapunov operator defined on $\mathcal{M}_{J+1}(\mathbb{C})$ by $\mathcal{L}_W(X)=WX+XW^T$, for all $X$. Thus, the discrete scheme leads to a Lyapunov algebraic recursive system. We now state the first result on the solvability of the numerical scheme.
\begin{Theorem}\label{theorem1}
	The system (\ref{Lyap1}) is uniquely solvable whenever $U^{0}$ and $U^{1}$ are known.
\end{Theorem}
In \cite{Benmabrouk-Ayadi2}, the authors have transformed the Lyapunov operator obtained from the discretization method into a standard linear operator acting on one column vector by juxtaposing the columns of the matrix $X$ horizontally which leads to an equivalent linear operator characterized by a fringe-tridiagonal matrix. We used standard computation to prove the invertibility of such an operator. Here. we do not apply the same computations as in \cite{Benmabrouk-Ayadi2}, but we develop different arguments. We will instead apply a differential calculus and topology technique (See \cite{HenriCartan} for example) to prove theorem \ref{theorem1}.
\begin{Lemma}\label{LemmeInversion}
Let $E$ be a finite dimensional ($\mathbb{R}$ or $\mathbb{C}$) vector space and $(\Phi_{n})_{n}$ be a sequence of endomorphisms converging uniformly to an invertible endomorphism $\Phi$. Then, there exists $n_{0}$ such that, for any $n\geq\,n_{0}$, the endomorphism $\Phi_{n}$ is invertible.
\end{Lemma}
Indeed, consider the set $Isom(E)$ of isomorphisms on $E$. It is regarded as the reciprocal image $det^{-1}(\mathbb{C}^*)$ with the determinant function. As this function is continuous, thus it consists of an open set in the set $L(E)$ of endomorphisms of $E$. Thus, as $\Phi\in\,Isom(E)$ there exists a ball $B(\Phi,r)\subset\,Isom(E)$. The elements $\Phi_n$ are in this ball for large values of $n$. So these are invertible.

Assume now that $l=o(h^{2+\varepsilon})$, with $\varepsilon>0$ which is always possible. Then, the coefficients appearing in $A$ and $W$ will satisfy as $h\longrightarrow0$ the following.
\begin{equation*}
A_n{j,j}=\displaystyle\frac{i}{2}-4\alpha_n\displaystyle\frac{l}{h^2}\longrightarrow\displaystyle\frac{i}{2};\quad\forall\,1\leq\,j\leq\,J.
\end{equation*}
For $1\leq\,j\leq\,J-1$,
\begin{equation*}
A_{j,j-1}=A_{j,j+1}=\displaystyle\frac{A_{0,1}}{2}=\displaystyle\frac{A_{J,J-1}}{2}=2\alpha_n\displaystyle\frac{l}{h^2}\longrightarrow0.
\end{equation*}
Next, observing that for all $X$ in the space $\mathcal{M}_{(J+1)^2}(\mathbb{C})$,
\begin{equation*}
\begin{matrix}
\|(\mathcal{L}_{A_n}-iI)(X)\| & =\|(A_n-\frac{i}{2}I)X+X(A_n^T-\frac{i}{2}I)\|\hfill
\cr &\leq2\|A-\frac{1}{2}I\|\|X\|,\hfill
\end{matrix}
\end{equation*}
it results that
\begin{equation}  \label{LWAtendsVersId}
\|\mathcal{L}_{A_n}-iI\|\leq\,Ch^{\varepsilon}.
\end{equation}
Consequently, the Lyapunov endomorphism $\mathcal{L}_{A_n}$ converges uniformly to the isomorphism $iI$ as $h$ goes towards 0 and $l=o(h^{2+\varepsilon})$ with $\varepsilon>0$. Using Lemma \ref{LemmeInversion}, the operator $\mathcal{L}_{A_n}$ is invertible for $h$ small enough. Hence, Theorem \ref{theorem1} is proved.
\section{Consistency, Stability and Convergence}
The consistency of the proposed method is done by evaluating the local truncation error arising from the discretization of the system
$$
iu_{t}+\Delta\,u+u-|u|^{-2\theta}u=0,\quad\in\Omega\times(t_0,+\infty)
$$
The principal part is
\begin{equation}\label{consistency1}
\mathcal{L}(x,y,t)=i\frac{l^{2}}{6} \frac{\partial^{3} u}{\partial t^{3}} + (\alpha - \gamma) l \Delta u_{t} + \frac{\alpha+\gamma}{2} l^2 \Delta u_{tt} + \beta \frac{h^{2}}{12}(\frac{\partial^{4} u}{\partial x^{4}} + \frac{\partial^{4} u}{\partial y^{4}} )
\end{equation}
 \begin{Lemma}
 \begin{itemize}
 \item $\alpha\neq\gamma$ and $l=o(h^2)$ the scheme is consistent with order $(h^2+l^2)$.
\item $\alpha=\gamma$ the scheme is unconditionally consistent with order $(h^2+l^2)$.
\end{itemize}
 \end{Lemma}
It is clear that the two operators $\mathcal{L}_{u}$ tend toward 0 as $l$ and $h$ tend to 0, which ensures the consistency of the method. Furthermore, the method may be always chosen to be consistent with an order 2 in time and space.

We now prove the stability of the method by applying the Lyapunov criterion which states that a linear system $\mathcal{L}(x_{n+1},x_{n},x_{n-1},\dots)=0$ is stable in the sense of Lyapunov if for any bounded initial solution $x_{0}$ the solution $x_{n}$ remains bounded for all $n\geq0$. Here, we will precisely prove the following result.
\begin{Lemma}\label{LyapunovStabilityLemma} $\mathcal{P}_{n}$: The solution $U^{n}$ is bounded independently of $n$ whenever the initial solution $U^{0}$ is bounded.
\end{Lemma}
We will proceed by induction on $n$. Assume firstly that $\|U^0\|\leq\eta$ for some $\eta$ positive. Using equation (\ref{Lyap1}), we obtain
\begin{equation}\label{LyapunovStability1}
\begin{matrix}
\|\mathcal{L}_{A_n}(U^{n+1})\|\leq\|\mathcal{L}_{B_n}\|.\|U^{n}\|+\|\mathcal{L}_{C_n}\|.\|U^{n-1}\|+\displaystyle\frac{\|F^{n-1}\|+\|F^n\|}{2}\hfill%
\end{matrix}%
\end{equation}
Next, recall that, for $l=o(h^{\varepsilon+2})$ small enough, $\varepsilon>0$, we have for $h\longrightarrow0$ the following uniform estimates (as the sequences $\alpha_n$, $\beta_n$ and $\gamma_n$ are bounded).
$$
\begin{array}{c}
B_n(j,j)=-2\sigma\beta_n\sim\,Ch^{\varepsilon}\longrightarrow0,\\
B_n(j,j-1)=B_n(j,j+1)=\sigma\beta_n\sim\,Ch^{\varepsilon}\longrightarrow0,\\
B_n(0,1)=B_n(J,J-1)=2\sigma\beta_n\sim\,Ch^{\varepsilon}\longrightarrow0.\\
\end{array}
$$
$$
\begin{array}{c}
C_n(j,j)=\psi_n\sim\,\displaystyle\frac{i}{2}+Ch^{\varepsilon}\longrightarrow\displaystyle\frac{i}{2},\\
C_n(j,j-1)=C_n(j,j+1)=\sigma\gamma_n\sim\,Ch^{\varepsilon}\longrightarrow0,\\
C_n(0,1)=C_n(J,J-1)=2\sigma\gamma_n\sim\,Ch^{\varepsilon}\longrightarrow0.\\
\end{array}
$$
As a consequence, for $h$ small enough,
\begin{equation}\label{LtildeBB}
\|\mathcal{L}_{B_n}\|\leq2\|B_n\|\leq\,Ch^{\varepsilon},
\end{equation}
and
\begin{equation}\label{LtildeBB}
\|\mathcal{L}_{C_n}\|\leq2\|C_n\|\leq\,\displaystyle\frac{1}{2}+Ch^{\varepsilon},
\end{equation}
We shall next use the following lemma deduced from (\ref{LWAtendsVersId}).
\begin{Lemma}\label{LWABounded}
For $h$ small enough, it holds for all $X\in \mathcal{M}_{(J+1)^{2}}(\mathbb{R})$ that
\begin{equation*}
\displaystyle\frac{1}{2}\|X\|\leq(1-Ch^{\varepsilon})\|X\|\leq\|\mathcal{L}_{A_n}(X)\|\leq(1+Ch^{\varepsilon})\|X\|\leq\displaystyle\frac{3}{2}\|X\|.
\end{equation*}
\end{Lemma}
Indeed, recall that equation (\ref{LWAtendsVersId}) affirms that $\|\mathcal{L}_{A_n}-iI\|\leq\,Ch^{\varepsilon}$ for some constant $C>0$. Consequently, for any $X$ we get
\begin{equation*}
(1-Ch^{\varepsilon})\|X\|\leq\|\mathcal{L}_{A_n}(X)\|\leq(1+Ch^{\varepsilon})\|X\|.
\end{equation*}
For $h\leq\displaystyle\frac{1}{C^{1/\varepsilon}}$, we obtain
\begin{equation*}
\displaystyle\frac{1}{2}\leq(1-Ch^{\varepsilon})<(1+Ch^{\varepsilon})\leq\displaystyle\frac{3}{2}
\end{equation*}
and thus Lemma \ref{LWABounded}. As a result, (\ref{LyapunovStability1}) yields that
\begin{equation}\label{LyapunovStability3-1}
\|U^{n+1}\|\leq\|(1+2Ch^{\varepsilon})\|U^{n}\|+2(1+Ch^{\varepsilon})\|U^{n-1}\|+\|U^{n-1}\|^{1-2\theta}+\|U^n\|^{1-2\theta}.
\end{equation}%
For $n=0$, this implies that
\begin{equation}
\|U^{1}\|\leq\|(1+2Ch^{\varepsilon})\|U^{0}\|+2(1+Ch^{\varepsilon})\|U^{-1}\|+\|U^{-1}\|^{1-2\theta}+\|U^0\|^{1-2\theta}.
\label{LyapunovStability3-2}
\end{equation}%
Using the discrete approximation
\begin{equation*}
U^{-1}=U^{0}-il(\Delta\,u_0+f(u_0))
\end{equation*}
and the fact that $u_0$ is sufficiently regular and thus bounded on the domain $\Omega$, we get
\begin{equation}
\|U^{-1}\|\leq\|U^{0}\|+Cl\leq\|U^{0}\|+Ch^{2+\varepsilon}.\label{U-1Bounds}
\end{equation}%
Hence, equation (\ref{LyapunovStability3-2}) yields that
\begin{equation}
\|U^{1}\|\leq(3+Ch^{\varepsilon})\|U^{0}\|+Ch^{2+\varepsilon}(1+Ch^{\varepsilon})+2(\|U^{0}\|+Ch^{2+\varepsilon})^{1-2\theta}.\label{LyapunovStability3-3}
\end{equation}%
Now, the Lyapunov criterion for stability states exactly that
\begin{equation}\label{LyapunovStability1}
\forall\,\,\varepsilon>0,\,\exists\,\eta>0\,\,\,s.t;\,\,\|U^{0}\|\leq\eta\,\,\Rightarrow\,\,\|U^{n}\|\leq\varepsilon,\,\,\forall\,n\geq0.
\end{equation}
For $n=1$ and $\|U^{1}\|\leq\varepsilon$, we seek an $\eta>0$ for which $\|U^{0}\|\leq\eta$. Indeed, using (\ref{LyapunovStability3-3}), this means that, it suffices to find $\eta $ such that
\begin{equation}
(3+Ch^{\varepsilon})\|U^{0}\|+Ch^{2+\varepsilon}(1+Ch^{\varepsilon})+2(\|U^{0}\|+Ch^{2+\varepsilon})^{1-2\theta}<\varepsilon.  \label{LyapunovStability3-4}
\end{equation}%
Choosing $h$ small enough ($h\leq\displaystyle\frac{\eta}{C})$, we seek eta such that
\begin{equation}
2\eta(\eta+2)+2^{2-2\theta}\eta^{1-2\theta}<\varepsilon.  \label{LyapunovStability3-5}
\end{equation}%
which is possible as the quantity at the left hand tends to 0 when $\eta\rightarrow0$.\\ Assume now that the $U^{k}$ is bounded for $k=1,2,\dots,n$ (by $\varepsilon_{1}$) whenever $U^{0}$ is bounded by $\eta$ and let $\varepsilon>0$. We shall prove that it is possible to choose $\eta$ satisfying $\|U^{n+1}\|\leq\varepsilon$. Indeed, from (\ref{LyapunovStability3-1}), we have
\begin{equation}\label{LyapunovStability3Ordren-1}
\|U^{n+1}\|\leq(3+Ch^{\varepsilon})\varepsilon_{1}+2\varepsilon_{1}^{1-2\theta}.
\end{equation}%
So, one seeks, $\varepsilon_{1}$ for which $(3+Ch^{\varepsilon})\varepsilon_{1}+2\varepsilon_{1}^{1-2\theta}<\varepsilon$ which is always possible. \\
Next, the convergence is a consequence of the well known Lax-Richtmyer equivalence theorem, which states that for consistent numerical approximations, stability and convergence are equivalent. Recall here that we have already proved in (\ref{consistency1}) that the used scheme is consistent. Next, Lemma \ref{LyapunovStabilityLemma}, Lemma \ref{LWABounded} and equation (\ref{LyapunovStability1}) yields the stability of the scheme. Consequently, the Lax equivalence Theorem guarantees the convergence. So as the following Lemma.
\begin{Lemma}\label{laxequivresult}
As the numerical scheme is consistent and stable, it is then convergent.
\end{Lemma}
\section{Numerical Implementation: 2-Particles Interaction}
In this section we propose to develop a numerical example to illustrate the efficiency of the numerical scheme proposed and studied previously. It consists of a model of interaction of two particles or two waves. We consider the inhomogeneous problem
\begin{equation}\label{SchrodingerNumericalExample}
\left\{\begin{matrix}\medskip i\displaystyle\frac{\partial u}{\partial
t}+\Delta\,u+u-|u|^{-2\theta}u=g(x,y,t)&\hbox{in}&\mathcal{Q},\hfill\cr\medskip
u(x,y,0)=v(x)v(y)&\hbox{in}&\overline{\Omega},\hfill\cr\medskip
\displaystyle\frac{\partial u}{\partial
n}(x,y,t)=0&\hbox{on}&\partial\Omega\times[0,T]\hfill\end{matrix}\right.
\end{equation}
where
$$
\Omega=]-L_0,L_0[^2,\ \mathcal{Q}=\Omega\times[0,T],\quad\,v(x)=\cos^2\left(\displaystyle\frac{\pi}{2L_0}x\right)
$$
and
$$
g(x,y,t)=\exp\left(-\displaystyle\frac{2i\pi^2}{L_0^2}t\right)\Bigl[\displaystyle\frac{\pi^2}{2L_0^2}(v(x)+v(y))+v(x)v(y)-v^{1-2\theta}(x)v^{1-2\theta}(y)\Bigr].
$$
The exact solution is
$$
u(x,y,t)=\exp\left(-\displaystyle\frac{2i\pi^2}{L_0^2}t\right)\cos^2\left(\displaystyle\frac{\pi}{2L_0}x\right)\cos^2\left(\displaystyle\frac{\pi}{2L_0}y\right).
$$

In the following tables, numerical results are provided. We computed for different space and time steps the discrete $L_2$-error estimates defined as follows. \begin{equation*}
\|X\|_{2}=\left( \displaystyle\sum_{i,j=1}^{N}|X_{ij}|^{2}\right)^{1/2}
\end{equation*}
for a matrix $X=(X_{ij})\in \mathcal{M}_{N+2}(\mathbb{C})$. Denote $u^{n}$ the net function $u(x,y,t^{n})$ and $U^{n}$ the numerical solution. The discrete $L_2$-error is
\begin{equation}
Er=\displaystyle\max_{n}\|U^{n}-u^{n}\|_{2}  \label{Er}
\end{equation}
on the grid $(x_{i},y_{j})$, $0\leq\,i,j\leq\,J+1$. We compute also the relative error between the exact solution and the numerical one as
\begin{equation}
Relative\,Er=\displaystyle\max_{n}\displaystyle\frac{\|U^{n}-u^{n}\|_{2}}{\|u^{n}\|_{2}}  \label{Errelative}
\end{equation}
on the same grid.

The domain $\Omega$ is chosen to be $\Omega=]-2\pi,2\pi[$. The time interval is $[0,1]$ for a choice $t_0=0$ and $T=1$. The following results are obtained for different values of $h$, $l$ and for $\theta=\displaystyle\frac{1}{4}$. We choed finally the barycenter calibrating parameters $\alpha_n$, $\beta_n$ and $\gamma_n$ to be
$$
\alpha_n=\displaystyle\frac{1}{4}+\displaystyle\frac{1}{2^{n+3}}\;\;\,\;\;\beta_n=\displaystyle\frac{1}{4}-\displaystyle\frac{1}{2^{n+3}}\;\;\hbox{and}\;\;
\gamma_n=\displaystyle\frac{1}{2}.
$$
\begin{table}[ht]
\begin{center}
\centerline{Table 1.}
\begin{tabular}{||l|l|l|l|l|l|l||}
\hline\hline
J  & $l$   & $\log(l)/\log(h)$ & $Er$ & $Relative\,Er$ & $Er/(l^2+h^2)$ \\ \hline
10 & 1/100 & -20.15 & $7,50.10^{-4}$ & 4,00 & $4,74.10^{-4}$\\ \hline
16 & 1/120 & 19.81  & $3,90.10^{-4}$ & 4,00 & $6,33.10^{-4}$\\ \hline
20 & 1/200 & 11.40  & $1,87.10^{-4}$ & 4,00 & $4,74.10^{-4}$\\ \hline
24 & 1/220 & 8.33   & $1,42.10^{-4}$ & 4,00 & $5,18.10^{-4}$\\ \hline
30 & 1/280 & 6.47   & $8,92.10^{-5}$ & 4,00 & $5,08.10^{-4}$\\ \hline
40 & 1/400 & 5.17   & $4,68.10^{-5}$ & 4,00 & $4,74.10^{-4}$\\ \hline
50 & 1/500 & 4.50   & $3,00.10^{-5}$ & 4,00 & $4,74.10^{-4}$\\ \hline\hline
\end{tabular}%
\end{center}
\end{table}
\newpage
As we see in the table, the numerical scheme converges with a convergence rate of $(l^2+h^2)$. Notice from the last column that the quantity $\displaystyle\frac{Er}{(l^2+h^2)}$ is of the order of $10^{-4}$ even in the case where the hypothesis $l=o(h^{2+\varepsilon})$ is not satisfied (line 1 of the table). In the 3rd column of the table, we notice the contribution of approximate solution to the exact solution. A relative error of $4\%$ meaning that the ratio $\displaystyle\frac{Unumerical}{Uexact}$ is of the order of $1\pm0.04$.
\begin{Remark}
The barycenter parameters $\alpha_n$, $\beta_n$ and $\gamma_n$ are applied to calibrates the position of the approximated solution relatively to the exact one. These parameters affect surely the numerical solution as well as the error estimates. In existing works such as \cite{Benmabrouk-Ayadi1}, \cite{Benmabrouk-Ayadi2}, \cite{Bratsos1}, \cite{Bratsos2}, \cite{Bratsos-et-al1} these are constants. Here. we adopted instead variable coefficients which may be generated by random procedures. These calibrations permits the use of implicit/explicit schemes by using suitable values. For example for $\alpha_n=\gamma_n=\frac{1}{2}$ and $\beta_n=0$, the barycentre estimation becomes
$$
U^{n}(\alpha_n,\beta_n,\gamma_n)=\displaystyle\frac{U^{n+1}+U^{n-1}}{2}
$$
which is an implicit estimation that guarantees an error of order $2$ in time.
\end{Remark}
\section{Conclusion}
This paper investigated the solution of the well-known NLS equation in two-dimensional case by applying a two-dimensional finite difference discretization. The continuous problem is firstly recasted into an algebraic discrete system involving Lyapunov-Syslvester matrix terms by using a full time-space discretization. Solvability, consistency, stability and convergence are then established by applying Lax-Richtmyer equivalence theorem and Lyapunov stability and by examining the Lyapunov-Sylvester operators. The method was finally improved by developing a numerical example issued from 2-particles interaction. It was shown to be efficient by means of error estimates.

\end{document}